\RequirePackage{ifpdf}
\ifpdf % We are running pdfTeX in pdf mode
\documentclass[pdftex]{sigma}
\else
\documentclass{sigma}
\fi

\begin{document}

\allowdisplaybreaks
\renewcommand{\PaperNumber}{002}

\FirstPageHeading

\renewcommand{\thefootnote}{$\star$}

\ShortArticleName{Raising and Lowering Operators for Askey--Wilson Polynomials}

\ArticleName{Raising and Lowering Operators\\ for Askey--Wilson Polynomials\footnote{This paper is a contribution 
to the Vadim Kuznetsov Memorial Issue ``Integrable Systems and Related Topics''.
The full collection is available at 
\href{http://www.emis.de/journals/SIGMA/kuznetsov.html}{http://www.emis.de/journals/SIGMA/kuznetsov.html}}}

\Author{Siddhartha SAHI}
\AuthorNameForHeading{S. Sahi}
\Address{Department of Mathematics, Rutgers University, New Brunswick, NJ 08903, USA} 
\Email{\href{mailto:sahi@math.rutgers.edu}{sahi@math.rutgers.edu}} 

\ArticleDates{Received September 20, 2006, in f\/inal form
December 27, 2006; Published online January 04, 2007}

\Abstract{In this paper we describe two pairs of raising/lowering 
operators for Askey--Wilson polynomials, 
which result from constructions involving very dif\/ferent techniques. 
The f\/irst technique is quite elementary, and depends only on the ``classical'' 
properties of these polynomials, viz.\ the $q$-dif\/ference equation 
and the three term recurrence. The second technique is less elementary, 
and involves the one-variable version of the double af\/f\/ine 
Hecke algebra.}


\Classification{33D45; 33D52; 33D80}

\begin{flushright}
\it Dedicated to the memory of Vadim Kuznetsov
\end{flushright}

\section{Introduction}

One of the approaches to proving integrality of the $(q,t)$-Kostka
coef\/f\/icients is the idea, due to Kirillov--Noumi \cite{kin,kin1} and
Lapointe--Vinet \cite{lv,lv1}, of using raising operators for Macdonald
polynomials. (See also \cite{gr,gt,kn,sa1} for other
approaches.) In their proof Kirillov and Noumi give an explicit construction
of such raising operators for the Macdonald polynomials $J_{\lambda}\left(
x;q,t\right)  $ for the root system of type $A_{n-1}$. They also pose the
problem of f\/inding analogous operators for the six-parameter Koornwinder
corresponding to the $BC_{n}$ root system.

This question was also raised by Tom Koornwinder at the Edinburgh conference
on symmetric functions organized by Vadim Kuznetsov. The case $n=1$
corresponds to the celebrated Askey--Wilson polynomials and Koornwinder's paper
\cite{ko1} from that conference contains partial results in this direction as
well as a survey of earlier results.

In this paper we such construct raising/lowering operators for Askey--Wilson
polynomials. In fact we describe \emph{two} such pairs of operators, which
result from constructions involving very dif\/ferent techniques.

The f\/irst technique is quite elementary, and depends only on the ``classical''
properties of these polynomials, \textit{viz.} the $q$-dif\/ference equation and
the three term recurrence. Therefore it can be applied to all the polynomials
in the Askey scheme. After this work was completed, we obtained a recent
preprint by T.~Koornwinder \cite{ko2}, the main result of which is very close
to this approach. Also through \cite{ko2} we discovered still earlier work
of G.~Bangerezako \cite{b} which obtains similar operators based on an
\textit{ad-hoc} factorization of the Askey--Wilson operator. Our approach
however is more direct and quite short.

The second technique is less elementary and involves the one-variable version
of the powerful Hecke algebra method as described in \cite{m,n,
ns,sa2,sa3,st}. This approach is related to a
fairly remarkable mathematical object~-- the double af\/f\/ine Hecke algebra
(see \cite{c1, c2,e,sa2}). The calculations, while
non-trivial to carry out, are conceptually rather straightforward. The
raising/lowering operators so obtained are dif\/ferent from those coming from
the ``classical'' method. This method also provides a new factorization of the
Askey--Wilson operator described in Lemma~\ref{factor}, which is much simpler
than that of Bangerezako.

In subsequent work, we hope to extend these methods to construct raising
operators for Koornwinder polynomials \cite{ko,v,sa2,sa3}.

\section{The classical approach}

\subsection{Askey Wilson polynomials}

The $q$-hypergeometric series is given by the formula%
\[
{}_{r}\phi_{s}\left(  \left.
\genfrac{}{}{0pt}{}{a_{1},\ldots,a_{r}}{b_{1},\ldots,b_{s}}%
\right|  q;y\right)  =\sum_{k\geq0}\frac{\left(  a_{1},\ldots,a_{r}\right)
_{k}}{\left(  b_{1},\ldots,b_{s}\right)  _{k}}\left(  -1\right)  ^{\left(
1+s-r\right)  k}q^{\left(  1+s-r\right)  \binom{k}{2}}\frac{y^{k}}{\left(
q\right)  _{k}},
\]
where the ``$q$-Pochhammer symbols'' are def\/ined by
\begin{gather*}
\left(  a,b,c,\ldots\right)  _{k}   :=\left(  a\right)  _{k}\left(
b\right)  _{k}\left(  c\right)  _{k}\cdots,\\
\left(  a\right)  _{k}    :=\left(  1-a\right)  \left(  1-aq\right)
\cdots\big(  1-aq^{k-1}\big).
\end{gather*}

The Askey--Wilson polynomials \cite{aw} are def\/ined by the formula%
\[
P_{n}\left(  z;a,b,c,d|q\right)  =\,\frac{\left(  ab,ac,ad\right)  _{n}}%
{a^{n}\left(  abcdq^{n-1}\right)  _{n}}\,_{4}\phi_{3}\left(  \left.
\genfrac{}{}{0pt}{}{q^{-n},abcdq^{n-1},az,az^{-1}}{ab,ac,ad}%
\right|  q;q\right).
\]
Since $\left(  q^{-n}\right)  _{k}$ vanishes for $k>n$, we have%
\[
_{4}\phi_{3}\left(  \cdots\right)  =\sum_{k=0}^{n}\left[  \frac{\left(
abcdq^{n-1}\right)  _{k}}{\left(  ab,ac,ad\right)  _{k}}\right]  \left[
\frac{\left(  q^{-n}\right)  _{k}q^{k}}{\left(  q\right)  _{k}}\right]
\left(  az\right)  _{k}\left(  az^{-1}\right)  _{k}.
\]

It follows that $P_{n}$ is a Laurent polynomial of degree $n,$ which is
moreover symmetric in $z$ and~$z^{-1}$ and is of the form
\[
P_{n}\left(  z;a,b,c,d|q\right)  =\left(  z^{n}+z^{-n}\right)  +\ \text{lower
terms}.
\]
It is also symmetric in $\left\{  a,b,c,d\right\}  $ although this is not
entirely obvious from the formula above.

We have chosen to normalize $P_{n}$ in order to make it monic. Of course there
are several other possible normalizations, and we discuss some of these below.

First of all, we remark that formula (3.1.7) of \cite{ks} considers the
polynomial
\[
\frac{\left(  ab,ac,ad\right)  _{n}}{a^{n}}\,_{4}\phi_{3}\left(  \left.
\genfrac{}{}{0pt}{}{q^{-n},abcdq^{n-1},az,az^{-1}}{ab,ac,ad}%
\right|  q;q\right)
\]
which is $\left(  abcdq^{n-1}\right)  _{n}$ times our $P_{n}.$

Next, since the Askey--Wilson polynomial is symmetric in $z$, $z^{-1},$ it can be
expressed as an (ordinary) polynomial of degree $n$ in%
\[
x=\left(  z+z^{-1}\right)  /2.
\]
The function $p_{n}\left(  x\right)  $ considered in (3.1.5) of \cite{ks}, is
monic in $x$, and hence is related to our normalization $P_{n}$ by the formula%
\[
p_{n}\left(  \frac{z+z^{-1}}{2}\right)  =2^{-n}P_{n}\left(  z\right).
\]

Finally, the polynomials $P_{n}$ are orthogonal with respect the inner product
$\left\langle \cdot ,\cdot \right\rangle $ def\/ined in (3.1.2) of \cite{ks}. If we def\/ine
\begin{gather}
Q_{n}    =\gamma_{n}P_{n},\label{gamma-n}
\end{gather}
where 
\begin{gather*}
\gamma_{n}    =\frac{(abq^{n},acq^{n},adq^{n},bcq^{n}
,bdq^{n},cdq^{n},q^{n+1})_{\infty}}{\left(  abcdq^{2n}\right)  _{\infty}
}\left(  abcdq^{n-1}\right)  _{n}.
\end{gather*}
Then $Q_{n}$ is dual to $P_{n}$ in the sense that
\[
\left\langle P_{m},Q_{n}\right\rangle =\delta_{m,n}.
\]

\subsection{Raising and lowering operators}

The main result of this section are the following raising and lowering
operators for the Askey--Wilson polynomials:

\begin{theorem}
\label{one} For all $n>1,$ the Askey--Wilson polynomials satisfy the relations%
\begin{gather*}
\left[  D\left(  z+z^{-1}\right)  -\lambda_{n-1}\left(  z+z^{-1}\right)
-\alpha_{n}\left(  \lambda_{n}-\lambda_{n-1}\right)  \right]  P_{n} 
=\left(  \lambda_{n+1}-\lambda_{n-1}\right)  P_{n+1},\\
\left[  D\left(  z+z^{-1}\right)  -\lambda_{n+1}\left(  z+z^{-1}\right)
-\alpha_{n}\left(  \lambda_{n}-\lambda_{n+1}\right)  \right]  Q_{n}  
=\left(  \lambda_{n-1}-\lambda_{n+1}\right)  Q_{n-1},
\end{gather*}
where $D$, $\lambda_{n}$ and $\alpha_{n}$ are as in \eqref{Ddef} and
\eqref{alpha-n} below.
\end{theorem}

\begin{proof}
The proof involves two key properties of the Askey--Wilson polynomials.

The f\/irst property is the `$q$-dif\/ference equation' from (3.1.7) ) of
\cite{ks} which asserts that $P_{n}$ is an eigenfunction for the Askey--Wilson
operator, i.e.
\begin{equation}
DP_{n}\left(  z\right)  =\lambda_{n}P_{n}\left(  z\right)  . \label{q-diff}
\end{equation}
The operator and its eigenvalue are def\/ined by
\begin{gather}
D    =A\left(  z\right)  \left(  T_{q}-1\right)  +A\left(  z^{-1}\right)
\left(  T_{q^{-1}}-1\right), \label{Ddef}\\
\lambda_{n}    =\left(  q^{-n}-1\right)  \left(  1-abcdq^{n-1}\right)
=\left(  q^{-n}+abcdq^{n-1}\right)  -\left(  1+abcdq^{-1}\right), \nonumber
\end{gather}
where $A\left(  z\right)  $ is the rational function%
\begin{equation}
A\left(  z\right)  =\frac{\left(  az,bz,cz,dz\right)  _{1}}{\left(
z^{2}\right)  _{2}}=\frac{\left(  1-az\right)  \left(  1-bz\right)  \left(
1-cz\right)  \left(  1-dz\right)  }{\left(  1-z^{2}\right)  \left(
1-qz^{2}\right)  } \label{Az}%
\end{equation}
and $T_{q}$ is the shift operator%
\[
T_{q}f\left(  z\right)  =f\left(  qz\right)  .
\]
(To forestall possible confusion we emphasize that, in accordance with custom,
we think of $f\left(  z\right)  $ as a Laurent polynomial rather than as a
function of $z$. This means that we have $T_{q}\left(  z^{k}\right)
=q^{k}z^{k}$ rather than $T_{q}\left(  z^{k}\right)  =q^{-k}z^{k}.$)

The second key property of these polynomials is the `normalized recurrence
relation' from (3.1.5) of \cite{ks} which can be rewritten in the form%
\begin{equation}
\left(  z+z^{-1}\right)  P_{n}=P_{n+1}+\alpha_{n}P_{n}+\frac{\gamma_{n-1}%
}{\gamma_{n}}P_{n-1}\qquad \text{for} \quad n>1, \label{recurr}%
\end{equation}
where
\begin{equation}
\alpha_{n}=a+1/a-\frac{a\left(  bcq^{n-1},bdq^{n-1},cdq^{n-1},q^{n}\right)
_{1}}{\left(  abcdq^{2n-2}\right)  _{2}}-\frac{\left(  abq^{n},acq^{n}
,adq^{n},abcdq^{n-1}\right)  _{1}}{a\left(  abcdq^{2n-1}\right)  _{2}}.
\label{alpha-n}%
\end{equation}
We combine these two properties as follows:

First apply the operators $D-\lambda_{n-1}$ and $D-\lambda_{n+1}$,
respectively, to the recurrence relation to get%
\begin{gather*}
\left(  D-\lambda_{n-1}\right)  \left(  z+z^{-1}-\alpha_{n}\right)  P_{n}  
=\left(  \lambda_{n+1}-\lambda_{n-1}\right)  P_{n+1},\\
\left(  D-\lambda_{n+1}\right)  \left(  z+z^{-1}-\alpha_{n}\right)  P_{n}  
=\frac{\gamma_{n-1}}{\gamma_{n}}\left(  \lambda_{n-1}-\lambda_{n+1}\right)
P_{n-1}.
\end{gather*}
Finally simplify, using the $q$-dif\/ference equation (\ref{q-diff}), to get%
\begin{gather*}
\left[  D\left(  z+z^{-1}\right)  -\lambda_{n-1}\left(  z+z^{-1}\right)
-\alpha_{n}\left(  \lambda_{n}-\lambda_{n-1}\right)  \right]  P_{n}  
=\left(  \lambda_{n+1}-\lambda_{n-1}\right)  P_{n+1},\\
\left[  D\left(  z+z^{-1}\right)  -\lambda_{n+1}\left(  z+z^{-1}\right)
-\alpha_{n}\left(  \lambda_{n}-\lambda_{n+1}\right)  \right]  Q_{n}  
=\left(  \lambda_{n-1}-\lambda_{n+1}\right)  Q_{n-1}
\end{gather*}
as desired.
\end{proof}

\section{The Hecke algebra approach}

In this section we provide raising/lowering operators for Askey--Wilson
polynomials based on Hecke algebra considerations \cite{sa2,sa3}. Once
again the main idea is quite straightforward, although the calculations are a
little more intricate. The resulting formulas are dif\/ferent and perhaps
slightly simpler.

\subsection{The Hecke algebra}

The key to this approach are the involutions $s_{1}$, $s_{0}$ which act on
Laurent polynomials as follows:%
\[
s_{1}f\left(  z\right)  =f\left(  z^{-1}\right) \qquad  \text{and}\qquad s_{0}f\left(
z\right)  =f\left(  qz^{-1}\right).
\]
Once again we regard these operators as acting on polynomials, rather than
functions, so that we have%
\[
s_{1}\left(  z^{k}\right)  =z^{-k}\qquad \text{and} \qquad s_{0}\left(  z^{k}\right)
=q^{k}z^{-k}.
\]
These operators provide a factorization of the $q$-shift operators, and one
has%
\[
s_{1}s_{0}=T_{q}\qquad \text{and} \qquad s_{0}s_{1}=T_{q^{-1}}.
\]

The af\/f\/ine Hecke algebra \cite{sa2,sa3} is the algebra of operators
generated by the two opera\-tors~$T_{0}$ and $T_{1}$ def\/ined as
\[
T_{i}:=t_{i}+r_{i}\left(  s_{i}-1\right),
\]
where
\begin{gather}
t_{0}    =-cd/q,\qquad r_{0}=\frac{\left(  z-c\right)  \left(  z-d\right)
}{\left(  z^{2}-q\right)  },\label{rrtt}\\
t_{1}    =-ab,\qquad r_{1}=\frac{\left(  1-az\right)  \left(  1-bz\right)
}{\left(  1-z^{2}\right)  }.\nonumber
\end{gather}

\begin{remark}
\label{Remark} The operator $T_{i}$ as def\/ined here is $t_{i}^{1/2}$ times the
corresponding operator from~\cite{sa2,sa3}. This accounts for the
slight dif\/ference between the formulas here and in~\cite{sa3}.
\end{remark}

From the def\/inition of $T_{1}$ it follows that a polynomial $f$ is symmetric
in $z$, $z^{-1}$, if and only if%
\begin{equation}
T_{1}f=t_{1}f. \label{symm}%
\end{equation}
Consequently, if $g$ is any polynomial, then the quadratic relation
(\ref{quad1}) implies that $\left(  T_{1}+1\right)  g$ is a symmetric
polynomial. The operators $T_{i}$ are deformations of $s_{i}$ and satisfy a
quadratic relation. For the convenience of the reader unfamiliar with the
Hecke algebra, we give a proof this relation.

\begin{lemma}
The operators $T_{i}$ satisfy the relation
\begin{equation}
\left(  T_{i}-t_{i}\right)  \left(  T_{i}+1\right)  =0. \label{quad1}%
\end{equation}
\end{lemma}

\begin{proof}
Def\/ine $s_{i}\left(  r_{i}\right)  =r_{i}^{\prime}$, then we claim that
\begin{equation}
r_{i}+r_{i}^{\prime}=t_{i}+1 \label{rrt1}%
\end{equation}
To see this, we calculate for $i=0,$%
\begin{gather*}
r_{0}+r_{0}^{\prime}    =\frac{\left(  z-c\right)  \left(  z-d\right)
}{\left(  z^{2}-q\right)  }+\frac{\left(  qz^{-1}-c\right)  \left(
qz^{-1}-d\right)  }{\left(  q^{2}z^{-2}-q\right)  }\\
\phantom{r_{0}+r_{0}^{\prime}}{}  =\frac{\left(  z-c\right)  \left(  z-d\right)  }{\left(  z^{2}-q\right)
}-\frac{\left(  q-cz\right)  \left(  q-dz\right)  }{q\left(  z^{2}-q\right)
}\\
\phantom{r_{0}+r_{0}^{\prime}}{}  =\frac{q\left(  z^{2}-cz-dz+cd\right)  -\left(  q^{2}-qcz-qdz+cdz^{2}%
\right)  }{q\left(  z^{2}-q\right)  }\\
\phantom{r_{0}+r_{0}^{\prime}}{}  =\frac{qz^{2}+qcd-q^{2}-cdz^{2}}{q\left(  z^{2}-q\right)  }=\frac{\left(
q-cd\right)  \left(  z^{2}-q\right)  }{q\left(  z^{2}-q\right)  }\\
\phantom{r_{0}+r_{0}^{\prime}}{}  =1-\frac{cd}{q}=1+t_{0}.
\end{gather*}
The calculation for $i=1$ is similar and simpler.

Now the quadratic relation can be proved as follows:%
\begin{gather*}
\left(  T_{1}-t_{1}\right)  \left(  T_{1}+1\right)     =r_{i}\left(
s_{i}-1\right)  \left[  t_{i}+1+r_{i}\left(  s_{i}-1\right)  \right]
\nonumber\\
\phantom{\left(  T_{1}-t_{1}\right)  \left(  T_{1}+1\right)}{}  =r_{i}\left[  \left(  t_{i}+1\right)  \left(  s_{i}-1\right)  +\left(
s_{i}r_{i}\right)  \left(  s_{i}-1\right)  -r_{i}\left(  s_{i}-1\right)
\right] \nonumber\\
\phantom{\left(  T_{1}-t_{1}\right)  \left(  T_{1}+1\right)}{}
  =r_{i}\left[  \left(  t_{i}+1\right)  \left(  s_{i}-1\right)  +\left(
r_{i}^{\prime}s_{i}\right)  \left(  s_{i}-1\right)  -r_{i}\left(
s_{i}-1\right)  \right] \nonumber\\
\phantom{\left(  T_{1}-t_{1}\right)  \left(  T_{1}+1\right)}{}
  =r_{i}\left[  \left(  t_{i}+1\right)  \left(  s_{i}-1\right)
+r_{i}^{\prime}\left(  1-s_{i}\right)  -r_{i}\left(  s_{i}-1\right)  \right]
\nonumber\\
\phantom{\left(  T_{1}-t_{1}\right)  \left(  T_{1}+1\right)}{}
  =r_{i}\left(  t_{i}+1-r_{i}-r_{i}^{\prime}\right)  \left(  s_{i}-1\right)
 \, \overset{{\rm by}\  \eqref{rrt1}}{=}\, 0. \tag*{\qed} % \tag{by \ref{rrt1}}%
\end{gather*}\renewcommand{\qed}{}
\end{proof}

The following result is an immediate consequence

\begin{corollary}
The operators $T_{i}$ are invertible, with%
\begin{equation}
t_{i}T_{i}^{-1}=T_{i}-t_{i}+1. \label{quad2}%
\end{equation}
\end{corollary}

We will also need a number of commutation results between the $T_{i}$ and the
multiplication operator by~$z.$ They follow directly from the def\/inition, and
we leave the (easy) proof to the reader.

\begin{lemma}
\label{commute} The operators\ $T_{i}$ satisfy the following commutation
relations%
\begin{gather*}
zt_{0}T_{0}^{-1}    =qT_{0}z^{-1}+c+d,\\
\left(  T_{1}+1\right)  z^{-1}    =t_{1}z^{-1}+zT_{1}+a+b,\\
\left(  T_{1}+1\right)  z    =z+t_{1}z^{-1}T_{1}^{-1}-\left(  a+b\right).
\end{gather*}
\end{lemma}

\subsection[Nonsymmetric Askey-Wilson polynomials]{Nonsymmetric Askey--Wilson polynomials}

The next ingredient in the Hecke algebra method are the nonsymmetric
Askey--Wilson polynomials. These are certain Laurent polynomials, $E_{n}$,
$n\in\mathbb{Z}$, which can be characterized up to multiples as eigenfunctions
of the operator%
\[
Y=T_{1}T_{0}.
\]
More precisely, one has%
\begin{gather}
YE_{n}    =\mu_{n}E_{n},\label{mun}
\end{gather}
where 
\begin{gather*}
\mu_{n}   =\left\{
\begin{array}{lll}
q^{n} & \text{for} & n<0,\\
q^{n}t_{1}t_{0}=q^{n-1}abcd & \text{for} & n\geq0.
\end{array}
\right. 
\end{gather*}

The symmetric Askey--Wilson polynomials $P_{\left|  n\right|  }$ are closely
related to $E_{\pm n}$. Up to normalization, one has $P_{0}=E_{0}=1$, while
for $n>0$ we have up to a scalar
\begin{equation}
P_{\left|  n\right|  }\sim\left(  T_{1}+1\right)  E_{\pm n}=c_{n}^{\pm}
E_{n}+c_{-n}^{\pm}E_{-n}. \label{PE}
\end{equation}
The explicit formula for the coef\/f\/icients $c_{n}^{\pm}$ and $c_{-n}^{\pm}$ is
known, but will not be needed in what follows.

We now def\/ine a slight variant of the Askey--Wilson operator, as follows:
\begin{equation}
D^{\prime}=A\left(  z\right)  \left(  T_{q}-s_{1}\right)  +A\left(
z^{-1}\right)  \left(  T_{q^{-1}}s_{1}-1\right).  \label{Dprime}
\end{equation}

Observe that $D$ and $D^{\prime}$ have the same action on functions which are
symmetric $z$ and $z^{-1},$ thus the Askey--Wilson polynomials satisfy%
\[
D^{\prime}P_{n}=DP_{n}=\lambda_{n}P_{n}.
\]

Just as the operator $s_{0}$ and $s_{1}$ factorize the $q$-shift operator, it
turns out that the opera\-tors~$T_{0}$ and $T_{1}$ provide a factorization of
$D^{\prime}$.

\begin{lemma}
\label{factor}The operator $D^{\prime}$ of formula \eqref{Dprime} admits the
following factorization:
\[
D^{\prime}=\left(  T_{1}+1\right)  \left(  T_{0}-t_{0}\right).
\]
\end{lemma}

\begin{proof}
To prove this, we calculate as follows
\begin{gather*}
\left(  T_{1}+1\right)  \left(  T_{0}-t_{0}\right)     =\left(  t_{1}
+r_{1}\left(  s_{1}-1\right)  +1\right)  \left(  r_{0}\left(  s_{0}-1\right)
\right) \,
\overset{{\rm by} \eqref{rrt1}}{=} \, \left(  r_{1}^{\prime}+r_{1}s_{1}\right)  \left(  r_{0}s_{0}-r_{0}\right)
\\
\phantom{\left(  T_{1}+1\right)  \left(  T_{0}-t_{0}\right)}{} 
 =r_{1}^{\prime}r_{0}\left(  s_{0}-1\right)  +r_{1}\widetilde{r}_{0}\left(
s_{1}s_{0}-s_{1}\right) 
  =r_{1}^{\prime}r_{0}\left(  T_{q}^{-1}s_{1}-1\right)  +r_{1}\widetilde
{r}_{0}\left(  T_{q}-s_{1}\right), \nonumber
\end{gather*}
where $r_{1}^{\prime}=r_{1}\left(  z^{-1}\right)  $ and $\widetilde{r}%
_{0}=r_{0}\left(  z^{-1}\right)  .$

Comparing the formulas for $r_{i}$ (\ref{rrtt}) and $A\left(  z\right)  $
(\ref{Az})$,$ we conclude that%
\[
A\left(  z\right)  =r_{1}\widetilde{r}_{0}\qquad {\rm and} \qquad A\left(  z^{-1}\right)
=r_{1}^{\prime}r_{0}%
\]
which completes the proof.
\end{proof}

\subsection{Raising and lowering operators}

To state our main result we need some notation. We write
\begin{equation}
e_{1}=a+b+c+d,\qquad e_{3}=abc+abd+acd+bcd. \label{e13}
\end{equation}

Also recall that for $n\geq0$, $\lambda_{n}$ is the symmetric Askey--Wilson
eigenvalue as in (\ref{Ddef}). For $n<0$ we def\/ine%
\[
\lambda_{n}=\lambda_{\left|  n\right|  }%
\]
and for all integral $n$ we set%
\begin{equation}
\beta_{n}=\frac{\lambda_{n}+1-\mu_{n-1}}{\mu_{n-1}-\mu_{-n}}e_{1}-\frac
{1-\mu_{n-1}}{\mu_{n-1}-\mu_{-n}}e_{3}. \label{beta-n}
\end{equation}

\begin{theorem}
\label{sraise} The Askey--Wilson polynomials satisfy the following relations:
\begin{gather}
\left[  D^{\prime}z+\left(  1-q^{1-n}\right)  \left(  z+z^{-1}\right)
+\beta_{-n}\right]  P_{n}    =\left(  q^{n}abcd-q^{1-n}\right)
P_{n+1},\qquad n\geq0,\label{sraise1}\\
\left[  D^{\prime}z+\left(  1-q^{n}abcd\right)  \left(  z+z^{-1}\right)
+\beta_{n}\right]  Q_{n}    =\left(  q^{1-n}-q^{n}abcd\right)  Q_{n-1},\qquad
n>0. \label{sraise2}%
\end{gather}
\end{theorem}

The key for the proof is the ``af\/f\/ine intertwiner'' for the nonsymmetric
Askey--Wilson polynomials from \cite{sa3}. This involves the additional
parameters $u_{0}$ and $u_{1},$ which satisfy the relations%
\[
a=t_{1}^{1/2}u_{1}^{1/2},\qquad b=-t_{1}^{1/2}u_{1}^{-1/2},\qquad c=q^{1/2}t_{0}^{1/2}%
u_{0}^{1/2},\qquad d=-q^{1/2}t_{0}^{-1/2}u_{0}^{1/2}.
\]

Now from Theorem 1.2 of \cite{sa3} we have, up to a multiple,%
\begin{gather}
E_{n}    \sim\left(  a_{n}U_{0}+b_{n}\right)  E_{-n-1},\label{eraise}
\end{gather}
where 
\[
a_{n}    =\big(  q^{\overline{n}}-q^{\overline{-n-1}}\big)
\qquad {\rm and} \qquad b_{n}    =q^{\overline{n}}\big(  u_{0}^{-1/2}-u_{0}^{1/2}\big)
+q^{-1/2}\big(  u_{1}^{-1/2}-u_{1}^{1/2}\big) 
\]
with 
\[
q^{\overline{n}}    =\mu_{n}t_{1}^{1/2}t_{0}^{1/2}
\]
and $U_{0}$ is the operator
\[
U_{0}=q^{-1/2}t_{0}^{1/2}T_{0}^{-1}z.
\]

We will derive Theorem \ref{sraise} from formula (\ref{eraise}); however some
remarks are in order before we proceed:

\begin{enumerate}\itemsep=0pt
\item  There is a typo in the statement of formula (\ref{eraise}) in Theorem
1.2 of \cite{sa3}, namely $n$ and $-n-1$ have been inadvertently switched.
This is easily seen by comparison with Theorem~4.1 from which Theorem~1.2 is derived.

\item  The formula for $U_{0}$ here is slightly dif\/ferent from that in
\cite{sa3} because of the dif\/ference in~$T_{0}$~-- see Remark~\ref{Remark}).

\item  Although Theorem 1.2 in \cite{sa3} is only stated (and needed) for
$n\geq0,$ it is easy to see that after the correction above it holds for all
integer $n$.

\item  Finally, we note that the ideas of \cite{sa2} and \cite{sa3} work in
the more general setting of Koornwinder polynomials, and they involve the
af\/f\/ine intertwiner $S_{0}$, which can also be written as
\[
S_{0}=\left[  Y,z^{-1}T_{1}^{-1}\right].
\]
It is expected that this operator will play a key role in the raising
operators for Koornwonder polynomials.

We now give the proof \ of Theorem \ref{sraise}.
\end{enumerate}

\begin{proof}
We f\/irst simplify (\ref{eraise}) by multiplying through by $q^{1/2}t_{1}%
^{1/2}t_{0}$. This gives
\begin{gather*}
E_{n}    \sim\left[  \left(  \mu_{n}-\mu_{-n-1}\right)  t_{0}T_{0}^{-1}%
z-\mu_{n}\left(  c+d\right)  -t_{0}\left(  a+b\right)  \right]  E_{-n-1}\\
\phantom{E_{n}}{} \sim\left[  t_{0}T_{0}^{-1}z-\frac{\mu_{n}\left(  c+d\right)  +t_{0}\left(
a+b\right)  }{\mu_{n}-\mu_{-n-1}}\right]  E_{-n-1}.
\end{gather*}
Replacing $n$ by $n-1,$ we get%
\[
E_{n-1}\sim\left[  t_{0}T_{0}^{-1}z-\frac{\mu_{n-1}\left(  c+d\right)
+t_{0}\left(  a+b\right)  }{\mu_{n-1}-\mu_{-n}}\right]  E_{-n}.
\]
For ease in subsequent calculations, we write this as
\begin{equation}
E_{n-1}\sim\left(  t_{0}T_{0}^{-1}z-\frak{\kappa}_{n}\right)  E_{-n},
\label{kapparaise}
\end{equation}
where
\begin{equation}
\frak{\kappa}_{n}=\frac{\mu_{n-1}y+t_{0}x}{\mu_{n-1}-\mu_{-n}} ,\qquad
x=a+b,\qquad y=c+d. \label{kappaxy}
\end{equation}

The key idea to obtain a raising operator for $P_{n}$ is as follows: By
formula (\ref{PE}), for $n>1,$ $P_{\left|  n\right|  }$~is a combination of
$E_{n}$ and $E_{-n}$. We f\/irst kill of\/f the $E_{n}$ component. This can be
accomplished by applying the operator~$Y-\mu_{n}$ to $P_{n}$. However
it is more convenient (and equivalent) to apply the operator%
\[
t_{1}t_{0}\left(  Y^{-1}-\mu_{n}^{-1}\right)  =t_{1}t_{0}T_{0}^{-1}T_{1}%
^{-1}-\frac{t_{1}t_{0}}{\mu_{n}}%
\]
For $n\neq0$ we have $\frac{t_{1}t_{0}}{\mu_{n}}=\mu_{-n}$. Thus since $P_{n}$
is symmetric, formula (\ref{symm}) implies that up to a~non-zero multiple, one
has%
\begin{equation}
\left(  t_{0}T_{0}^{-1}-\mu_{-n}\right)  P_{n}\sim E_{-n}. \label{project}%
\end{equation}
Although the argument given above only applies for $n\neq0,$ it is easy to see
that formula (\ref{project}) is true (up to a non-zero multiple) for $n=0$ as
well! Now combining formulas (\ref{PE}), (\ref{kapparaise}), and
(\ref{project}) we conclude that up to a multiple, we have%
\[
P_{\left|  n-1\right|  }    \sim R_{n}P_{\left|  n\right|  },
\]
where 
\[
R_{n}    =\left(  T_{1}+1\right)  \left(  t_{0}T_{0}%
^{-1}z-\frak{\kappa}_{n}\right)  \left(  t_{0}T_{0}^{-1}-\mu_{-n}\right).
\]

The main problem now is to simplify the expression of the operator $R_{n}$
using properties of~$P_{\left|  n\right|  }$.

We f\/irst calculate using Lemma \ref{commute} and (\ref{quad2}), as follows%
\begin{gather*}
  \left(  t_{0}T_{0}^{-1}z-\frak{\kappa}_{n}\right)  \left(  t_{0}T_{0}%
^{-1}-\mu_{-n}\right) \\
\qquad{}  =t_{0}T_{0}^{-1}zt_{0}T_{0}^{-1}-t_{0}T_{0}^{-1}\left(  \mu_{-n}%
z+\frak{\kappa}_{n}\right)  +\kappa_{n}\mu_{-n}\\
\qquad{}  =t_{0}T_{0}^{-1}\left(  qT_{0}z^{-1}+y\right)  -t_{0}T_{0}^{-1}\left(
\mu_{-n}z+\frak{\kappa}_{n}\right)  +\kappa_{n}\mu_{-n}\\
\qquad{}  =t_{0}T_{0}^{-1}\left(  y-\mu_{-n}z-\frak{\kappa}_{n}\right)  +qt_{0}%
z^{-1}+\kappa_{n}\mu_{-n}\\
\qquad{}  =\left(  T_{0}-t_{0}+1\right)  \left(  y-\mu_{-n}z-\frak{\kappa}%
_{n}\right)  +\left(  qt_{0}z^{-1}+\kappa_{n}\mu_{-n}\right).
\end{gather*}
Applying $T_{1}+1$ to this, we get by Lemma \ref{factor}
\[
R_{n}=\left(  D^{\prime}+T_{1}+1\right)  \left(  -\mu_{-n}z+y-\frak{\kappa
}_{n}\right)  +\left(  T_{1}+1\right)  \left(  qt_{0}z^{-1}+\kappa_{n}\mu
_{-n}\right).
\]

To simplify this further we note that the commutation relations of Lemma
\ref{commute} can be rewritten as and hence imply%
\begin{gather*}
\left(  T_{1}+1\right)  z^{-1}    =t_{1}z^{-1}+zT_{1}+x,\\
\left(  T_{1}+1\right)  z    =z+t_{1}z^{-1}T_{1}^{-1}-x.
\end{gather*}
Also on $P_{\left|  n\right|  }$, $D^{\prime}$ acts by $\lambda_{n}$ while
$T_{1}$ acts by $t_{1}$. Therefore $R_{n}$ acts by the operator%
\begin{gather*}
  -\mu_{-n}\left[  D^{\prime}z+z+z^{-1}-x\right]  +\left(  \lambda_{n}%
+t_{1}+1\right)  \left(  y-\frak{\kappa}_{n}\right) \\
\qquad{} +qt_{0}\left[  t_{1}\left(  z+z^{-1}\right)  +x\right]  +\left(
t_{1}+1\right)  \kappa_{n}\mu_{-n}.
\end{gather*}

Dividing by $-\mu_{-n}$, we see that up to a multiple
\begin{equation}
\left[  D^{\prime}z+\left(  1-\frac{qt_{1}t_{0}}{\mu_{-n}}\right)  \left(
z+z^{-1}\right)  +\beta_{n}^{\prime}\right]  P_{\left|  n\right|  }\sim
P_{\left|  n-1\right|  }, \label{raise-lower}
\end{equation}
where
\[
\beta_{n}^{\prime}=-x-\frac{qt_{0}}{\mu_{-n}}x-\frac{\lambda_{n}+t_{1}+1}
{\mu_{-n}}\left(  y-\frak{\kappa}_{n}\right)  -\left(  t_{1}+1\right)
\frak{\kappa}_{n}.
\]

We now show that $\beta_{n}^{\prime}$ is equal to $\beta_{n}$ from formula
(\ref{beta-n}). For this we simplify the expression, substituting for
$\frak{\kappa}_{n},$ using (\ref{kapparaise}) above to get%
\begin{gather*}
\beta_{n}^{\prime}    =-x-\frac{qt_{0}}{\mu_{-n}}x+\frac{\left(  \lambda
_{n}+t_{1}+1\right)  }{\mu_{-n}}\frac{\mu_{-n}y+t_{0}x}{\mu_{n-1}-\mu_{-n}%
}-\left(  t_{1}+1\right)  \frac{\mu_{n-1}y+t_{0}x}{\mu_{n-1}-\mu_{-n}}\\
\phantom{\beta_{n}^{\prime}}{} =\frac{\lambda_{n}+1-\mu_{n-1}}{\mu_{n-1}-\mu_{-n}}y+\frac{1-\mu_{n-1}}%
{\mu_{n-1}-\mu_{-n}}t_{1}y+c_{1}x+c_{2}\left(  qt_{0}x\right),
\end{gather*}
where $c_{1}$ and $c_{2}$ are moderately complicated expressions which can be
computed explicitly. However, we can save some computation by observing that
since the result is \emph{a priori} symmetric in $\left\{  a,b,c,d\right\}  $,
$c_{1}$ and $c_{2}$ must reduce to the coef\/f\/icients of $y$ and $-t_{1}y$ respectively.

It follows then that we have%
\begin{gather*}
\beta_{n}^{\prime}    =\frac{\lambda_{n}+1-\mu_{n-1}}{\mu_{n-1}-\mu_{-n}
}\left(  x+y\right)  +\frac{1-\mu_{n-1}}{\mu_{n-1}-\mu_{-n}}\left(
qt_{0}x+t_{1}y\right) \\
\phantom{\beta_{n}^{\prime}}{}  =\frac{\lambda_{n}+1-\mu_{n-1}}{\mu_{n-1}-\mu_{-n}}e_{1}-\frac{1-\mu_{n-1}%
}{\mu_{n-1}-\mu_{-n}}e_{3}=\beta_{n}.
\end{gather*}

Replacing $n$ by $-n$, formula (\ref{raise-lower}) becomes
\[
\left[  D^{\prime}z+\left(  1-\frac{qt_{1}t_{0}}{\mu_{n}}\right)  \left(
z+z^{-1}\right)  +\beta_{-n}\right]  P_{\left|  n\right|  }\sim P_{\left|
-n-1\right|  }.
\]
For $n\geq0,$ we have $\mu_{n}=q^{n}t_{1}t_{0}$ and this becomes
\begin{equation}
\left[  D^{\prime}z+\left(  1-q^{1-n}\right)  \left(  z+z^{-1}\right)
+\beta_{-n}\right]  P_{n}\sim P_{n+1} \label{mraise1}%
\end{equation}
which is (\ref{sraise1}) up to a multiple.

For $n\geq1,$ we have $\mu_{-n}=q^{-n}$ and formula (\ref{raise-lower})
becomes%
\begin{equation}
\left[  D^{\prime}z+\left(  1-q^{n}abcd\right)  \left(  z+z^{-1}\right)
+\beta_{n}\right]  P_{n}\sim P_{n-1} \label{mraise2}%
\end{equation}
which is (\ref{sraise2}) up to a multiple.

It remains then only to calculate the multiples in (\ref{mraise1}),
(\ref{mraise2}).

To determine the multiple (\ref{mraise1}), it suf\/f\/ices to calculate the
coef\/f\/icient of $z^{n+1}$ on the left. For this we divide the left side of
(\ref{mraise1}) by $z^{n+1}$ and take the limit as $z\rightarrow\infty$. This
gives%
\begin{gather*}
  \lim_{z\rightarrow\infty}\frac{1}{z^{n+1}}\left(  A\left(  z\right)
q^{n+1}z^{n+1}-A\left(  z^{-1}\right)  z^{n+1}+\left(  1-q^{1-n}\right)
z^{n+1}\right) \\
\qquad{}  =\frac{abcd}{q}q^{n+1}-1+1-q^{1-n}=q^{n}abcd-q^{1-n}
\end{gather*}
which proves formula (\ref{sraise1}).

To determine the multiple in\ (\ref{mraise2}), we rewrite (\ref{mraise2}) in
the form%
\[
\left[  D^{\prime}z+\left(  1-q^{n}abcd\right)  \left(  z+z^{-1}\right)
+\beta_{n}\right]  P_{n}=c_{n}\left(  q^{1-n}-q^{n}abcd\right)  \frac
{\gamma_{n-1}}{\gamma_{n}}P_{n-1},\qquad n>0
\]
for some unknown constant $c_{n}$. It then suf\/f\/ices to show that $c_{n}=1.$
Subtracting this from (\ref{sraise1}) we have%
\begin{gather*}
  \left[  \left(  q^{n}abcd-q^{1-n}\right)  \left(  z+z^{-1}\right)
+\beta_{-n}-\beta_{n}\right]  P_{n}\\
\qquad  =\left(  q^{n}abcd-q^{1-n}\right)  \left(  P_{n+1}+c_{n}\frac{\gamma_{n-1}
}{\gamma_{n}}P_{n-1}\right)
\end{gather*}
or
\[
\left(  \text{ }z+z^{-1}\right)  P_{n}=P_{n+1}+\frac{\beta_{n}-\beta_{-n}
}{q^{n}abcd-q^{1-n}}P_{n}+c_{n}\frac{\gamma_{n-1}}{\gamma_{n}}P_{n-1}.
\]

Comparing this with the recurrence relation (\ref{recurr}) we deduce $c_{n}
=1$, which proves formula (\ref{sraise2}) and completes the proof of the theorem.

We note in passing that comparison with (\ref{recurr}) also proves the
following identity (which can be verif\/ied independently):%
\begin{gather*}
\alpha_{n}=\frac{\beta_{n}-\beta_{-n}}{q^{n}abcd-q^{1-n}}.\tag*{\qed}
\end{gather*}\renewcommand{\qed}{}
\end{proof}

\subsection*{Acknowledgements}

We would like to thank the (anonymous) referee for several insightful
suggestions which have improved the paper considerably. The referee has also
pointed out that one can give an alternative proof of formulas (\ref{sraise1})
and (\ref{sraise2}) by combining Theorem \ref{one} with the following identity
relating the operators $D$ and $D^{\prime}:$
\begin{gather*}
  \left[  (1-q^{2})D^{\prime}z+q^{2}D(z+z^{-1})-q(z+z^{-1})D\right]  f\\
\qquad  =(1-q)\left[  (e_{1}-e_{3})-(1-abcd)(z+z^{-1})\right]  f,
\end{gather*}
which holds for all symmetric Laurent polynomials $f$.

\pdfbookmark[1]{References}{ref}
\LastPageEnding

\begin{thebibliography}{99} 

\footnotesize\itemsep=-0.2pt

\bibitem{aw}Askey R., Wilson J., Some basic hypergeometric polynomials that generalize Jacobi polynomials, 
{\it Mem. Amer. Math. Soc.}
\textbf{319} (1985), 1--53.

\bibitem{b}Bangerezako G., The factorization method for the Askey--Wilson
polynomials, {\it J. Comput. Appl. Math.} \textbf{107} (1999), 219--232.

\bibitem{c1}Cherednik I., Double af\/f\/ine Hecke algebras,
Knizhnik--Zamolodchikov equations, and Macdonald's operators, {\it Int. Math.
Res. Not.} (1992), no.~9, 171--180.

\bibitem{c2}Cherednik I., Double af\/f\/ine Hecke algebras and Macdonald's
conjectures, {\it Ann. of Math.} \textbf{141} (1995), 191--216.

\bibitem{e}Etingof P., Oblomkov A., Rains E., Generalized double
af\/f\/ine Hecke algebras of rank 1 and quantized del Pezzo surfaces, 
\href{http://arxiv.org/abs/math.QA/0406480}{math.QA/0406480}.

\bibitem{gr}Garsia A., Remmel R., Plethystic formulas and
positivity for $q,t$-Kostka coef\/f\/icients, in Mathematical Essays in Honor of
Gian-Carlo Rota, Editors B.~Sagan and R. Stanley, {\it Progr. Math.} \textbf{161}
(1998), 245--262.

\bibitem{gt}Garsia A., Tesler G., Plethystic formulas for the
Macdonald $q,t$-Kostka coef\/f\/icients, {\it Adv. Math.} \textbf{123} (1996), 144--222.

\bibitem{is}Ion B.,  Sahi S., Triple groups and Cherednik algebras,
{\it Contemp. Math.} \textbf{417} (2006), 183--206, \href{http://arxiv.org/abs/math.QA/0304186}{math.QA/0304186}.

\bibitem{kin}Kirillov A., Noumi M., $q$-dif\/ference raising operators
for Macdonald polynomials and the integrality of transition coef\/f\/icients, in
Algebraic Methods and $q$-Special Functions, {\it CRM Proceedings and Lecture Notes}
{\bf 22} (1999), 227--243, \href{http://arxiv.org/abs/q-alg/9605005}{q-alg/9605005}.

\bibitem{kin1}Kirillov A., Noumi M., Af\/f\/ine Hecke algebras and
raising operators for Macdonald polynomials, {\it Duke Math. J.} \textbf{93}
(1998), 1--39,  \href{http://arxiv.org/abs/q-alg/9605004}{q-alg/9605004}.

\bibitem{kn}Knop F., Integrality of two variable Kostka functions, {\it J.
Reine Angew. Math.} \textbf{482} (1997), 177--189, \href{http://arxiv.org/abs/q-alg/9603027}{\mbox{q-alg/9603027}}.

\bibitem{ks}Koekoek R., Swarttouw R., The Askey-scheme of
hypergeometric orthogonal polynomials and its $q$-analogue, Delft University of
Technology, Department of Technical Mathematics and Informatics, Report no.
98-17 (1998), \url{http://aw.twi.tudelft.nl/~koekoek/askey/ch3/par1/par1.html}.

\bibitem{ko}Koornwinder T., Askey--Wilson polynomials for root systems
of type BC, {\it Contemp. Math.} \textbf{138} (1992), 189--204.

\bibitem{ko1}Koornwinder T., Lowering and raising operators for some
special orthogonal polynomials, {\it Contemp. Math.} \textbf{417} (2006), 227--238.

\bibitem{ko2}Koornwinder T., The structure relation for Askey--Wilson
polynomials, {\it J. Comput. Appl. Math.}, to appear, \href{http://arxiv.org/abs/math.CA/0601303}{math.CA/0601303}.

\bibitem{lv}Lapointe L., Vinet L., Creation operators for the
Macdonald and Jack polynomials, {\it Lett. Math. Phys.} \textbf{40} (1997), 269--286.

\bibitem{lv1}Lapointe L., Vinet L., Rodrigues formulas for the
Macdonald polynomials, {\it Adv. Math.} \textbf{130} (1997), 261--279, \href{http://arxiv.org/abs/q-alg/9607025}{q-alg/9607025}.

\bibitem{m}Macdonald I., Af\/f\/ine Hecke algebras and orthogonal polynomials,
Cambridge University Press, 2003.

\bibitem{n}Noumi M., Macdonald--Koornwinder polynomials and af\/f\/ine Hecke
algebras, {\it RIMS Kokyuroku} \textbf{919} (1995), 44--55 (in Japanese).

\bibitem{ns}Noumi M., Stokman J., Askey--Wilson polynomials: an
af\/f\/ine Hecke algebra approach, in Laredo Lectures on Orthogonal Polynomials
and Special Functions, Editors R.~Alvarez-Nodarse, F.~Marcellan and W.~Van Assche, 
 Nova Science Publishers, 2004, 111--144, \href{http://arxiv.org/abs/math.QA/0001033}{math.QA/0001033}.

\bibitem{sa1}Sahi S., Interpolation, integrality, and a generalization
of Macdonald's polynomials, {\it Int. Math. Res. Not.} (1996), no.~10, 457--471.

\bibitem{sa2}Sahi S., Nonsymmetric Koornwinder polynomials and
duality, {\it Ann. of Math.} \textbf{150} (1999),  267--282, \href{http://arxiv.org/abs/q-alg/9710032}{\mbox{q-alg/9710032}}.

\bibitem{sa3}Sahi S., Some properties of Koornwinder polynomials,
{\it Contemp. Math.} \textbf{254} (2000), 395--411.

%(South Hadley, MA, 1998), pp. 395-411, Contemp. Math. \textbf{254}, Amer. Math. Soc., Providence, RI, 2000

\bibitem{st}Stokman J., Koornwinder polynomials and af\/f\/ine Hecke
algebras, {\it Int. Math. Res. Not.} (2000), no.~19, 1005--1042, \href{http://arxiv.org/abs/math.QA/0002090}{math.QA/0002090}.

\bibitem{v}van Diejen J., Self-dual Koornwinder--Macdonald polynomials,
{\it Invent. Math.} \textbf{126} (1996), 319--339, \href{http://arxiv.org/abs/q-alg/9507033}{\mbox{q-alg/9507033}}.
\end{thebibliography}
\end{document}